\documentclass[12pt]{amsart}
\usepackage{amssymb}

\newtheorem{theorem}{Theorem}

\newtheorem{proposition}{Proposition}
\newtheorem{corollary}{Corollary}
\newtheorem{remark}{Remark}
\newtheorem{definition}{Definition}
\newtheorem{lemma}{Lemma}

\newcommand{\Z}{{\mathbb Z}}
\newcommand{\R}{{\mathbb R}}
\newcommand{\C}{{\mathbb C}}
\newcommand{\CP}{{\mathbb CP}}

\begin{document}

\author{V.A.~Vassiliev}
\address{Steklov Mathematical Institute of Russian Academy of Sciences;
\newline \hspace*{3mm}
National Research University Higher School of Economics} \email{vva@mi.ras.ru}

\thanks{Research supported by the Russian Science Foundation grant, project
16-11-10316}

\title{Multiplicities of bifurcation sets of Pham singularities}
\maketitle

\section{Introduction} Let $f:(\C^n,0) \to (\C,0)$ be a holomorphic function defined in a neighborhood of the orign in $\C^n$ and having an isolated singularity at $0$; let $F:(\C^n\times \C^k,0) \to (\C,0)$ be an arbitrary deformation of $f$, that is, a family of functions $f_\lambda$ depending holomorphically on the parameter $\lambda \in \C^k$ with $f_0\equiv f$, see \cite{AVGL}. 

\begin{definition}[see e.g. \cite{AVGL}] \rm  The {\em caustic} $\Delta \subset \C^k$ of the deformation $F$ is the set of all parameter values $\lambda \in \C^k$ such that the corresponding function $f_\lambda$ has a non-Morse critical point close to the origin in $\C^n$. The {\em Maxwell set}  of this deformation is the closure of the set of all parameter values $\lambda$ such that $f_\lambda$ has equal critical values at different critical points
close to the origin.
\end{definition}

\begin{definition} \rm The {\em mixed Stokes' set} (respectively, the {\em pure Stokes' set}) of deformation $F$ is the closure of the set of all parameter values $\lambda$ such that $f_\lambda$ has three different critical values $\alpha, \alpha', \beta$ satisfying the equality $$\alpha+\alpha'=2\beta ,$$ 
(respectively, has four different critical values $\alpha, \alpha', \beta, \beta'$ satisfying \begin{equation}\alpha+\alpha'=\beta+\beta'   \ \ ).\label{csdef}\end{equation} 
\end{definition}

If our deformation $F$ is large enough, that is, it is a {\em
versal} deformation of $f$, then the local multiplicities of the caustic, of the Maxwell set, and of both Stokes' sets
at the point $0 \in \C^k$ do not depend on the choice of this deformation. 
\medskip

We calculate the local multiplicities of the Maxwell sets and mixed Stokes' sets for all {\em Pham singularities}, that is, for singularities of the form \begin{equation} f(z_1, \dots, z_n) = z_1^{a_1+1}
+ \dots + z_n^{a_n+1}. \label{pham} \end{equation} Also we calculate the multiplicities of pure Stokes' sets for some Pham singularities, including all singularities in one variable. 

The multiplicities of caustics of Pham singularities are calculated in \cite{VFA} (see formula (\ref{caus}) below), where also a two-sided estimate of multiplicities of Maxwell sets was given. Our present result replaces this estimate by an equality, proving that its upper side is sharp. 

These multiplicities provide upper bounds for the complexity of certain programs enumerating all topologically distinct morsifications of complicated real function singularities. Moreover they allow us to write the stop rules in these programs, see \S \ref{motivat} at the end of this article. 

All these multiplicities depend semicontinuously on the singularity $f$, therefore our calculations give also upper estimates of them for arbitrary isolated holomorphic function singularities: the multiplicity of any of these kinds for a singularity $f$ does not exceed that for any Pham singularity such that $f$ occurs as its perturbation.

\section{Statements of main theorems}

\subsection{Maxwell set and caustic}

Assume that $a_1 \geq a_2 \geq \dots \geq a_n$ in (\ref{pham}), and define the number
\begin{equation}
L(a_1, \dots, a_n) \equiv \sum_{i=1}^n (a_1 \cdot \ \dots \ \cdot  a_{i-1})
(a_i^2-1) (a_{i+1} \cdot \ \dots \ \cdot a_n)^2. \label{Mset}
\end{equation}

\begin{theorem}
\label{ma2} The local multiplicities $C(f)$ and $M(f)$ of the caustic and the
Maxwell set of any versal deformation of the isolated singularity $f$ given by $($\ref{pham}$)$ satisfy the equation
\begin{equation}
3C(f) + 2M(f) = L(a_1, \dots, a_n).
\label{mnest}\end{equation}
\end{theorem}

On the other hand, according to \cite{VFA}, we have 
\begin{equation}
 C(f) = (a_1 \cdot \ \dots \ \cdot a_n)\left(\frac{a_1-1}{a_1}+ \dots +
 \frac{a_n-1}{a_n}\right) \ ,\label{caus}
\end{equation}
which allows us to express $M(f)$ from (\ref{mnest}). \medskip

{\bf A problem}. Is it correct that for an arbitrary isolated function singularity the local multiplicity of the caustic does not exceed $n\mu(f)$, where $\mu(f)$ is the Milnor number of $f$? 

\subsection{Mixed Stokes' set}

\begin{theorem} \label{mss}
The multiplicity at $0 \in \C^k$ of the mixed Stokes' set of the Pham singularity $($\ref{pham}$)$ is equal to 
\begin{equation}
\label{bbbb}
\binom{{a_1 \dots a_n}}{1,2}\frac{a_1+1}{a_1} + \sum_{j=2}^n (a_1 \dots a_{j-1})\binom{{a_j \dots a_n}}{1,2}\frac{a_{j-1}-a_j}{a_{j-1}a_j} \ ,
\end{equation}
where $\binom{A}{1,2}\equiv A(A-1)(A-2)/2$.
\end{theorem}

\subsection{Pure Stokes' set}

\begin{theorem}
\label{cs11}
If $n=1,$  $f=z^{d+1}$, and
$d$ is odd, then the multiplicity of the pure Stokes' set is equal to 
\begin{equation}
(d+1)(d-1)(d-2)(d-3)/8.
\end{equation}
\end{theorem}

\begin{theorem}
\label{cs12}
 If $n=1$, $f=z^{d+1},$ and $d$ is even, then the multiplicity of the pure Stokes' set is equal to 
\begin{equation}
\label{cs12f}
(d-2)((d+1)(d-1)(d-3)+1)/8.
\end{equation}
\end{theorem}

\begin{theorem} \label{cs2}Suppose that $n=2$ and  $f=x^{a+1}+y^{b+1},$ $a\geq b$, where both $a$ and $b$ are odd. Then the local multiplicity of the pure Stokes' set is equal to 
\begin{equation}
\label{anscM2}
\frac{a+1}{2a} \binom{ab}{{2,2}} + \frac{a-b}{2b}  \binom{b}{{2,2}}  + \frac{a-b}{a} \binom{a}{2} \binom{b}{2}(b-1)  + 
\frac{1}{a} \binom{a}{2} \binom{b}{2}. 
\end{equation}
\end{theorem}

{\bf Problem.} Calculate the multiplicities of all these sets in the terms of Newton diagrams for the generic functions $f$ with these diagrams.

\subsection{Corollaries for the  homogeneous case}

In the asymptotic formulas of the following proposition we assume that $n$ is fixed and the degrees are growing; recall that the Milnor number $\mu(f)$ of a Pham singularity (\ref{pham}) is equal to $a_1 \cdots a_n$.

\begin{proposition}
If $f$ has the form $($\ref{pham}$)$ 
and all exponents $a_1, \dots, a_n$ are equal to one another $($and are denoted by $a)$, then we have
$$C(f) = a^{n-1}(a-1) \sim n \mu(f),$$
$$M(f) = a^{n-1}((a+1)(a^n-1)-3n(a-1))/2 \sim \mu^2(f)/2;$$
the multiplicity of the mixed Stokes' set is $ \sim \mu^3(f)/2;$ \\
if additionally $n=1$, or $n=2$ and $a_1=a_2,$ then the multiplicity of the pure Stokes' set is $ \sim \mu^4(f)/8.$  
\end{proposition}

\section{Proof of Theorem \ref{ma2}.} 
\label{proofma2}

First of all, we replace the function (\ref{pham}) by
\begin{equation}
\label{pham1} f = {\small \frac{1}{a_1+1}} z_1^{a_1+1} + \dots +  {\small \frac{1}{a_n+1}} z_n^{a_n+1}
\end{equation}
(which can be done by the dilation of coordinates) because it simplifies the calculations very much.
We will use the versal deformation of this function (\ref{pham1}) consisting of all polynomial functions 
\begin{equation} \label{vers} f_\lambda\equiv f+ \sum_\alpha \lambda_\alpha z^\alpha,\end{equation} 
where $\alpha=(\alpha_1, \dots, \alpha_n) \in \Z^n$ are multi-indices with integer values 
\begin{equation}\alpha_i \in [0, a_i-1],\label{parall} \end{equation} 
$z^\alpha \equiv z_1^{\alpha_1} \cdot \ \dots \ \cdot z_n^{\alpha_n}$, and $\lambda_\alpha$ are the parameters, $\lambda = \{\lambda_\alpha\}$. In this case $k = \mu(f)=a_1 \cdot \ \dots \ \cdot a_n$.

Any function $f_\lambda$, where $\lambda \in \C^k \setminus \Delta$ is sufficiently close to $0\in \C^k$, has exactly $\mu(f)$ different critical points close to the origin in $\C^n$. 
Consider the complex-valued function on $\C^k\setminus \Delta$, whose value at any point $\lambda$ is equal to the product (over all
$2\binom{\mu(f)}{2}$ ordered pairs of different critical points of $f_\lambda$ close to the origin) of differences of the corresponding two critical values. This function is holomorphic and regular close to generic points of the caustic (these points $\lambda$ correspond to the functions $f_\lambda$ having exactly one non-Morse critical point of type $A_2$, while the non-generic points of the caustic form a set of complex codimension $\geq 2$ in $\C^k$). Therefore this function it can be extended to a holomorphic function $D$ in an entire neighborhood of the origin in $\C^k$. 

\begin{theorem} \label{mp} The degree of the restriction of the function $D$ to a generic line through the origin in $\C^k$
 is equal to $($\ref{Mset}$)$.
\end{theorem}

Theorem \ref{ma2} follows immediately from this one, because this function $D$ vanishes with multiplicity 3 at generic points of the caustic, and with multiplicity 2 at generic points of the Maxwell set.  \medskip

\subsection{Proof of Theorem \ref{mp}} 
\label{pro6}

Any line through the origin in $\C^k$ consists of functions of the form $f-\varepsilon \varphi$, where $\varepsilon \in \C$ is the parameter of the line, and $\varphi$ is a polynomial containing only the monomials $z^\alpha$ with $\alpha=(\alpha_1, \dots, \alpha_n)$ satisfying the conditions (\ref{parall}). We can and will assume that $\varphi(0)=0$, because $D(\lambda)=D(\lambda')$ if $\lambda$ and $\lambda'$ differ only in the coordinate $\lambda_0$, that is, $f_\lambda-f_{\lambda'}$ is a constant function. Let  \begin{equation}
\varphi_0=q_1z_1 + \dots + q_nz_n \label{linprt} \end{equation} 
be the linear part of $\varphi$. 

Consider the space $\CP^{k-2}$ of all lines of this form \begin{equation}\{f-\varepsilon \varphi\} \label{lines}\end{equation} in the subspace
$\C^{k-1} \subset \C^k$ distinguished by the last condition $\lambda_0=0$. The lines of this form, which are generic with respect to the function $D$ (that is, the degree at 0 of the restriction of $D$ to these lines is the minimal possible),
fill in a Zariski open subset in ${\mathbb CP}^{k-2}$. It is enough to prove the assertion of Theorem (\ref{mp}) for an arbitrary line from this subset, therefore we can and will assume that the linear part (\ref{linprt}) of our function $\varphi$ satisfies the following conditions:   
\begin{itemize}
\item $q_1=1$ (which can be made by the choice of the parameter $\varepsilon$);
\item all coefficients $q_i$ in (\ref{linprt}) are positive and $\leq 1$; and
\item $q_{i+1} \ll q_i$ if $a_{i+1} =a_i$.
\end{itemize}

Let us fix an arbitrary such function $\varphi$ and the corresponding line (\ref{lines}).
It is not obvious apriori that the line consisting of functions $f-\varepsilon \varphi_0$, where $\varphi_0$ is the linear part of $\varphi$, is also generic with respect to the function $D$. (Moreover, in \S \ref{proCS} we will see that in a similar problem concerning the pure Stokes' sets no perturbation $f-\varepsilon \varphi$ with linear $\varphi$ can be generic). Still we will calculate the degree of the restriction of the function $D$ to such a line, and find that this degree is equal to (\ref{Mset}). Then we connect $\varphi_0$ and $\varphi$ by a sequence of polynomials $\varphi_1, \dots, \varphi_{n-1}, \varphi_n \equiv \varphi$ and prove that this degree is the same for any two consecutive terms of this sequence.  Namely, we define $\varphi_j$ as the sum of all monomials of the polynomial $\varphi$, which either are of degree 1, or depend on variables $z_1, \dots, z_j$ only.

\begin{lemma}
\label{mainll} 
For any $j=0,1, \dots, n$, the critical points of the function $f-\varepsilon \varphi_j$ with sufficiently small $|\varepsilon|$ can be split into $a_1$ collections ``of depth 1'' with $a_2  \cdots a_n$ points in each, and the distances between the critical values at these points from different collections decreasing as $\asymp |\varepsilon|^{(a_1+1)/a_1}$ when $\varepsilon$ tends to 0; any of these $a_1$ collections can be subdivided into $a_2$ collections ``of depth 2'' with $a_3 \cdots a_n$ critical points in each and distances between the critical values at the points from different such subcollections $($inside one collection of depth 1$)$ decreasing as $\asymp |\varepsilon|^{(a_2+1)/a_2} ,$ etc: for any $i=1, \dots, n-1$ any of $a_1 \cdots a_{i-1}$ collections of depth $i$ distinguished in the previous steps can be split into $a_i$ collections of depth $i+1$ with $a_{i+1}\cdots a_n$ critical points in each and distances between the values at the points from different such collections of depth $i+1$ $($inside one collection of depth $i)$ decreasing as $\asymp |\varepsilon|^{(a_i+1)/a_i}.$
\end{lemma}

In what follws it will be convenient to consider the set of all $a_1 \cdots a_n$ critical points of such a function as the collection ``of depth $0$''.

\unitlength 0.5mm \linethickness{0.4pt}
\begin{figure}
\begin{picture}(130,120)
\put(120,58){\circle*{0.7}}

\put(128,66){\circle*{2}} \put(128,50){\circle*{2}} \put(112,66){\circle*{2}}
\put(112,50){\circle*{2}}

\put(0,58){\circle*{0.7}}

\put(8,66){\circle*{2}} \put(8,50){\circle*{2}} \put(-8,66){\circle*{2}}
\put(-8,50){\circle*{2}}

\put(90,108){\circle*{0.7}}

\put(98,116){\circle*{2}} \put(98,100){\circle*{2}} \put(82,116){\circle*{2}}
\put(82,100){\circle*{2}}

\put(90,8){\circle*{0.7}}

\put(98,16){\circle*{2}} \put(98,0){\circle*{2}} \put(82,16){\circle*{2}}
\put(82,0){\circle*{2}}

\put(30,108){\circle*{0.7}}

\put(38,116){\circle*{2}} \put(38,100){\circle*{2}} \put(22,116){\circle*{2}}
\put(22,100){\circle*{2}}

\put(30,8){\circle*{0.7}}

\put(38,16){\circle*{2}} \put(38,0){\circle*{2}} \put(22,16){\circle*{2}}
\put(22,0){\circle*{2}}

\end{picture}
\caption{Critical values for $a_1=6, a_2=4$}
\label{critvalpic} 

\end{figure}
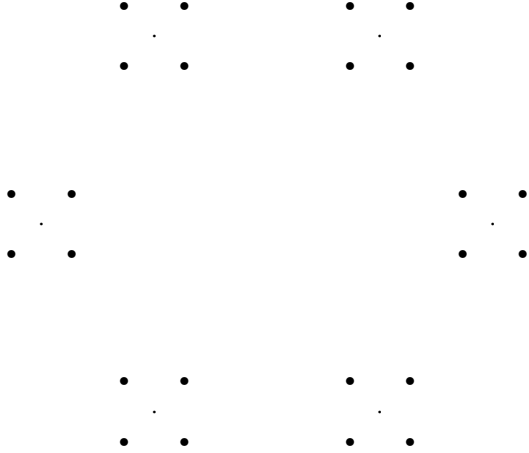

\begin{corollary} For any $j=0,1, \dots, n$, the degree of the restriction of the function $D$ to the line $\{f-\varepsilon \varphi_j\} \subset \C^{k-1}$ is equal to $($\ref{Mset}$)$.
\end{corollary}

{\it Proof.} For any $i=1, \dots, n$, the product of differences of critical values over all ordered pairs of critical points of $f-\varepsilon \varphi_j$, which belong to the same collections of depth $i-1$, but not of depth $i$, 
vanishes as $|\varepsilon|^{s_i}$, where $s_i$ is the $i$-th summand in (\ref{Mset}).  \hfill $\Box$ \medskip 

For $j=n$ the assertion of this corollary gives us Theorem \ref{mp}. \medskip

{\bf Basis of induction: proof of Lemma \ref{mainll} for $\varphi\equiv\varphi_0$}

The function $f-\varepsilon \varphi_0$ has $a_1\dots a_n$ critical points with coordinates 
\begin{equation} \label{valf0} z_1 = ( q_1\varepsilon) ^{1/a_1}, \dots,z_n =  (q_n\varepsilon)^{1/a_n}, \end{equation}
where any of expressions $( q_i\varepsilon)^{1/a_i}$ runs over all $a_i$ values of this root. 
The critical values at these points are equal to 
$$ \varepsilon \sum_{i=1}^n c_i \cdot (q_i\varepsilon)^{1/a_i},$$
where $c_i = -a_1/(a_1+1)q_i.$ Then the collections of depth $i$ of critical points from Lemma \ref{mainll} are just the points (\ref{valf0}) with coinciding coordinates $z_1, \dots, z_i$. \hfill $\Box$ \medskip

The step of induction in the proof of Lemma \ref{mainll} will be as follows: we show that the passage from $f-\varepsilon \varphi_{j-1}$ to $f-\varepsilon \varphi_j$ translates all collections of critical values of depth $j$ (i.e. consisting of $a_{j+1}  \cdots a_n$ elements) as whole bodies by the distances of size $O_{|\varepsilon| \to +0} (|\varepsilon|^{1+1/{a_j}+1/{a_1}}),$ thus not changing the orders of distances between the points from different subcollections of this depth (these orders are $\asymp (|\varepsilon|^{1+1/a_i})$, $i \leq j$) and not changing at all the distances between the critical values from one and the same such collection.

\begin{lemma}
\label{mainlll}
For $|\varepsilon|>0$  small enough, there is a natural and depending continuously on $\varepsilon$ one-to one correspondence between the sets of critical points of all functions $f-\varepsilon \varphi_j,$ $j=0, \dots, n$. For any $j, j' \in \{0, \dots, n\}$ 
the difference of values of any coordinate $z_i$, $i=1, \dots, n$, at the critical points of functions 
$f-\varepsilon \varphi_j$ and $f-\varepsilon \varphi_{j'}$ related by this correspondence decreases as $O_{|\varepsilon| \to +0}(|\varepsilon|^{1/{a_i}+1/{a_1}})$ .  
\end{lemma}
\medskip

{\it Proof.}
Consider the function $z_i^{a_i} - q_i$ in one variable, that is, the expression of the function $\partial(f-\varphi_0)/\partial z_i$. Its derivative is different from 0 at all its zero points $z_{0,i}$, i.e. at $a_i$th roots of $q_i$. Therefore any such point admits a boundary of some radius $T>0$ in $\C^1$, and a constant $c>0$ such that $|z_i^{a_i}-q_i|\geq c|z_i-z_{0,i}|$ if $|z_i-z_{0,i}|\leq T$. Let us choose these constants $T$ and $c$ in such a way that this inequality will hold for all $i=1, \dots, n$. Then by the dilation of the coordinate $z_i$ we obtain that 
\begin{equation}
\label{est2}
|z_i^{a_i}-q_i\varepsilon| \geq c |\varepsilon|^{\frac{a_i-1}{a_i}} | z_i-z_{0,i,\varepsilon}| \quad \mbox{if} \quad |z_i-z_{0,i,\varepsilon}|\leq |\varepsilon|^{\frac{1}{a_i}}T \ ,
\end{equation}    
where $z_{0,i,\varepsilon}$ is an arbitrary root of the polynomial $z_i^{a_i}-q_i\varepsilon$.

Let us fix some number $M>0$. 
For any critical point $Z_\varepsilon=(z_{0,1,\varepsilon}, \dots, z_{0,n,\varepsilon})$ of the function $f - \varepsilon \varphi_0$, consider its neighbourhood $U_M(Z_\varepsilon)$ containing all the points $(z_1, \dots, z_n)$ such that 
\begin{equation} \label{est3} |z_i -z_{0,i,\varepsilon}| \leq M |\varepsilon|^{1/a_i+1/a_1}\end{equation} 
for any $i=1, \dots, n$. If $\varepsilon$ is small enough then all these neighbourhoods 
belong to the polydisc where $|z_i|\leq 1$ for any $i$, and do not have common points, in particular any of them contains exactly one critical point (\ref{valf0}).  The boundary of this neighbourhood $U_M(Z_\varepsilon)$ consists of $n$ pieces, on any of which one of inequalities (\ref{est3}) becomes an equality.

Any function $\partial(\varphi_j-\varphi_0)/\partial z_k$ can be represented in the form $\sum_{i=1}^n z_i \psi_{i,j,k}(z)$ where $\psi_{i,j,k}$ are some polynomials. Let $C=\max_{i,j,k}\psi_{i,j,k}(z)$ where $i,j,k \in \{1, \dots, n\}$, $i \leq j$, and $z$ belongs to the polydisc where all $|z_i|\leq 1$. 
By the triangle inequality in any of our neighbourhoods $U_M(Z_\varepsilon)$ we have $|z_i| \leq |\varepsilon|^{1/a_i}+M|\varepsilon|^{1/a_i + 1/a_1},$ and hence

\begin{equation}
\label{est1}
\left|\frac{\partial(\varphi_j-\varphi_0)}{\partial z_i}\right| < nC(1+M|\varepsilon|^{1/a_i})|\varepsilon|^{1/a_1}.
\end{equation}

Suppose that \begin{equation} \label{estt}
c M > nC(1+M|\varepsilon|^{1/a_i}).\end{equation} 
Applying (\ref{est2}), (\ref{est3}) and (\ref{est1}) at the points of such a piece of the boundary of our neighborhood $U_M(Z_\varepsilon)$, where the $i$-th inequality (\ref{est3})  becomes an equality, we get 
\begin{equation}
\left|\frac{\partial (f-\varepsilon \varphi_0)}{\partial z_i}\right| \geq c M |\varepsilon|^{1+1/a_1} > nC(1+M|\varepsilon|^{1/a_i}) |\varepsilon|^{1+1/a_1} \geq \left| \frac{\partial(\varepsilon \varphi_j-\varepsilon \varphi_0)}{\partial z_i} \right|.
\label{est4}
\end{equation}
Therefore by the ``argument principle'' the indices of the vector fields $\mbox{grad} (f-\varepsilon \varphi_0)$ and $\mbox{grad} (f-\varepsilon \varphi_j)$ on the entire boundary of any such neighborhood coincide (and are equal to 1, as for the first of them).  In particular, there is exactly one critical point of the function $f-\varepsilon \varphi_j$ inside this neighborhood.

It remains to prove that all our assumptions used to get this conclusion are compatible, namely, the following statement. 

\begin{lemma} 
There exists $M>0$ such that for all $\varepsilon$ with $|\varepsilon|$ small enough the following three conditions are satisfied: the inequality $($\ref{estt}$)$, the right-hand condition in $($\ref{est2}$)$, and the condition that all our neighborhoods $U_M(Z_\varepsilon)$ of the critical points of $f-\varepsilon \varphi_0$ belong to the polydisc $\max |z_i|\leq 1$. 
\end{lemma}

{\it Proof.} We can choose $M=2nC/c$, and then impose three restrictions on $|\varepsilon|$: $nC|\varepsilon|^{1/a_i}< c/2$, $|\varepsilon|^{1/a_1}<T/M$ and $|\varepsilon|^{1/a_1}+M|\varepsilon|^{2/a_1}<1.$ These three conditions imply our three restrictions, respectively. 
\hfill $\Box$  $\Box$

\subsection{Induction step in the proof of Lemma \ref{mainll}}

Let $Z_\varepsilon(j-1)$ and $Z_\varepsilon(j)$ be some critical points of the functions 
$f-\varepsilon \varphi_{j-1}$ and $f-\varepsilon \varphi_j$, related to one another by the correspondence from Lemma \ref{mainlll}, that is, lying in one and the same neighbourhood $U_M(Z_\varepsilon)$ considered in the proof of this Lemma. We need to estimate the difference of the corresponding critical values, i.e. the number
$$ |\varepsilon||\varphi_j(Z_\varepsilon(j)) -\varphi_{j-1}(Z_\varepsilon(j-1))| . $$ Consider the family of functions 
\begin{equation}
\label{link}
f-\varepsilon(\varphi_{j-1}+t (\varphi_j-\varphi_{j-1})
\end{equation}
 depending on the parameter $t \in [0,1]$. By Lemma \ref{mainlll} any of these functions has a unique critical point in the domain $U_M(Z_\varepsilon)$. Define the function $u:[0,1]\to \C$, whose value at the point $t$ is the critical value of the corresponding function (\ref{link}) at this its critical points. The desired difference is equal to 
\begin{equation}
\label{int}
|u(1)-u(0)| = \left|\int_{0}^1 \frac{\partial u(t)}{\partial t} dt \right|. 
\end{equation}

\begin{lemma}[see e.g. Lemma 9.7.5 in \cite{LL}] The derivative in the 
integral $($\ref{int}$)$ is equal to the value of the function $-\varepsilon(\varphi_j -\varphi_{j-1})$ at the considered critical point of the function $f-\varepsilon(\varphi_{j-1}+t (\varphi_j-\varphi_{j-1}))$. \hfill $\Box$
\end{lemma} 

The absolute value of the function $\varphi_j-\varphi_{j-1}$ in $U_M(Z_\varepsilon)$  is uniformly estimated from above  
by $\mbox{ const} \times |\varepsilon|^{{1/a_1}+{1/a_j}}$, hence the difference (\ref{int}) decreases as $O_{\varepsilon \to 0}|\varepsilon|^{1+1/a_1+1/a_j}$.
Therefore moving from $\varphi_{j-1}$ to $\varphi_j$ does not affect the order of decrease of differences between the critical values from different collections of depth $\leq j$ (this order is equal to $|\varepsilon|^{1+1/a_i}$, where $i \in \{1, \dots, j\}$ is the first level of depth of collections separating these critical values).

On the other hand, this moving does not affect also the distances inside the collections of critical values of depth $m>j$. Indeed, any such collection of critical values of  the function $f-\varepsilon f_0$ is nothing else than the set of critical values of the function 
\begin{equation}\sum_{r=m+1}^n \left(\frac{1}{a_r+1} z_r^{a_r+1} - \varepsilon z_r\right)
\label{sumtail}
\end{equation} 
added to some single critical value of the function   $$\sum_{r=1}^m \left(\frac{1}{a_r+1} z_r^{a_r+1} - \varepsilon z_r\right).$$ The function $\varphi_j - \varphi_0$ does not contain any monomials depending on variables $z_r,$ $r \geq m$, therefore the corresponding collection 
for the function $f-\varepsilon \varphi_j$ is just the same set of critical values of the function (\ref{sumtail}) added to some critical value of the function in $m-1$ variables which is the sum of all monomials of the function $f-\varepsilon \varphi_j$ depending on these coordinates. Lemma \ref{mainll} is proved. \hfill $\Box$ $\Box$ $\Box$

\section{Proof of Theorem \ref{mss}} Consider again the parameter space $\C^k$ of the deformation (\ref{vers}) of $f$, remove from it its caustic $\Delta$, and define the following complex-valued function $Y$ on the remaining space.

Given a point $\lambda \in \C^k \setminus \Delta$, take all $\binom{\mu(f)}{1,2}$ choices of three critical points $Z_1, Z_2, Z_3$ of $f_\lambda$ (the first of which is distinguished, and the other two are not, so that the choice $(Z_1,Z_2, Z_3)$ is equal to $(Z_1, Z_3, Z_2)$), and define $Y(\lambda)$ as the product of all corresponding numbers \begin{equation}
\label{jjj}
2f_\lambda(Z_1) -f_\lambda(Z_2)-f_\lambda(Z_3).
\end{equation}

The obtained function can be extended to a holomorphic function close to the generic points of the caustic, and hence also to the whole neidhborhood of the origin in $\C^k$. It vanishes exactly on the mixed Stokes' set, and has multiplicity 1 at its generic points. Therefore the multiplicity of the mixed Stokes' set is equal to the degree of this function at the point $\lambda=0$ or, which is the same, of its restriction to a generic line through the origin in $\C^k$, consisting of polynomials $f -\varepsilon \varphi$, $\varepsilon \in \C^1$. The restriction of the function $Y$ to this line is a holomorphic function in the coordinate $\varepsilon$; let us calculate the order of its zero at the point $\varepsilon=0$. 

 Suppose first that $\varphi$ is linear and satisfies the conditons from \S \ref{pro6}.

The set of critical values of $f-\varepsilon \varphi$ has the structure described in Lemma \ref{mainll}, see also Fig.~\ref{critvalpic}. Consider these values as functions of $\varepsilon$. 
It is easy to see that the absolute value of the difference (\ref{jjj}) decreases as $|\varepsilon|^{1+1/a_j}$ with $\varepsilon \to 0$,  where $j$ is the smallest depth of collections of critical points such that not all three points $Z_1, Z_2, Z_3$ belong to one and the same collection of this depth. In particular, all $\binom{a_1 \cdots a_n}{1,2}$ factors of the function $Y(\varepsilon)$ give the contribution at least $1+1/a_1$ to the exponent of this function in $\varepsilon$.
These contributions give us the first summand in (\ref{bbbb}).
 Moreover, there are some $a_1 \times \binom{a_2\cdots a_n}{1,2}$ choices which give the contribution at least $1+1/a_2$: it are all possible triples which belong to one and the same collection of depth 1.
Accounting this addition (equal to $\frac{1}{a_2}-\frac{1}{a_1}=\frac{a_1-a_2}{a_1a_2}$) over all of these choices, we get the second summand in (\ref{bbbb}) (corresponding to $j=2$), etc. Finally, we get that the restriction of the function $Y$ to the generic line $\{f-\varepsilon \varphi\}$ with linear $\varphi$ is a monomial with the exponent equal to (\ref{bbbb}). 

Moreover, the same considerations as in \S \ref{proofma2} prove that replacing the generic linear function $\varphi$ by a non-linear one with the same linear part does not change the orders of decrease of all our factors (\ref{jjj}), so that the order (\ref{bbbb}) remains the same. \hfill $\Box$    
 
\begin{remark} \rm By analogy with the formula (\ref{bbbb}), the formula (\ref{Mset}) can be rewritten as 
$$2\binom{a_1 \cdots a_n}{2} \frac{a_1+1}{a_1} + 2 \sum_{j=2}^n a_1 \cdots a_{j-1} \binom{a_j \cdots a_n}{2} \frac{a_{j-1}-a_j}{a_{j-1}a_j} \ .$$
\end{remark}

\section{Proofs for pure Stokes' sets}

\label{proCS}

We proceed exactly as in the previous two sections and define the following function $\Omega(\lambda)$ on the space $\C^k \setminus \Delta$. Given $\lambda \in \C^k \setminus \Delta$, consider all possible choices of four unordered critical points $Z_1, \dots, Z_4 \in \C^n$ of $f_\lambda$ separated somehow into two pairs with the fixed order of these pairs: in total $\binom{\mu(f)}{2,2}$ choices. For any such choice we take the corresponding difference
\begin{equation}
\label{uuuu}
f_\lambda(Z_1) + f_\lambda(Z_2) - f_\lambda (Z_3)-f_\lambda(Z_4)
\end{equation} 
(i.e. the sum of critical values at the critical points from the first pair minus the similar sum for the second one) and define the function $\Omega(\lambda)$ as the product of these differences over all our choices. 
This function vanishes exactly on the pure Stokes' set with multiplicity 2. 

To calculate its degree, consider again a line through the origin in $\C^k$ consisting of functions $f-\varepsilon \varphi$, where $\varphi$ is a linear combination of monomials $z^\alpha$ with exponents $\alpha$ as in (\ref{vers}) and $\varphi(0)=0$. 

\subsection{Proof of theorem \ref{cs11}} 
Suppose that $n=1$ and $a_1 \equiv d$ is odd. Take first $\varphi\equiv z$. In this case all the factors (\ref{uuuu}) of the restriction of the function $\Omega$ to our line vanish exactly as $\varepsilon^{1+1/d}$. Therefore this restriction has a root of multiplicity $\frac{d+1}{d} \binom{d}{2,2} \equiv (d+1)(d-1)(d-2)(d-3)/4$. 

Moreover, the same arguments as in the proof of Theorem \ref{mp} show that in this case adding the non-linear monomials to $\varphi$ does not change the order of decrease of our factors and keeps the multiplicity of this root unchanged. The multiplicity of the pure Stokes' set is a half of this number, which proves Theorem \ref{cs11}. \hfill $\Box$

\subsection{Proof of theorem \ref{cs12}}

In the case of $n=1$ and even $d>2$, the line consisting of functions $f-\varepsilon \varphi$ with linear $\varphi$ 
(that is, of functions $z^{d+1}-\varepsilon z$)
will not be generic with respect to the pure Stokes' set. Indeed, if $(Z_1 , Z_2)$ and $(Z_3,Z_4)$ are two different pairs of opposite critical points of such a function then the factor (\ref{uuuu}) is equal to zero for all values of $\varepsilon$. Therefore let us choose the tentative function $\varphi$ in the form $z+\alpha z^2+\dots$, $\alpha \neq 0$.  

It is easy to calculate that in this case a majority of factors (\ref{uuuu})
still vanish as $\varepsilon^{1+1/a_1}$, and the remaining $\frac{d}{2} \left( \frac{d}{2}-1\right)$
many (corresponding to all possible choices of ordered pairs of unordered pairs of opposite roots of the polynomial $f'-\varepsilon$) vanish as $\varepsilon^{1+2/d}.$  Adding all these exponents and dividing the sum by 2 (i.e. by the multiplicity of the function $\Omega$ at the generic point of the pure Stokes' set) we get the number (\ref{cs12f}). \hfill $\Box$

\subsection{Proof of Theorem \ref{cs2}.} 

Now we have $f= x^{a+1} + y^{b+1}$, $a$ and $b$ odd, and $\varphi_0=x+q y$, where $q$ is positive, and $q \ll 1$ if $a=b$. 
The critical values of $f-\varepsilon \varphi_0$ split into $a$ collections, with $b$ points in each (see Fig. \ref{critvalpic}, where however the case of {\em even} $a$ and $b$ is shown).

We have $\binom{ab}{2,2}$ possible choices of quadruples of critical points $Z_1,\dots, $ $Z_4$ divided somehow into two numbered pairs. For a majority of these choices the expression (\ref{uuuu}) decreases as $\varepsilon^{1+1/a}$. However there are the following three possible degenerations increasing this exponent:

1) all points $Z_i$ belong to one and the same collection of $b$ values. This happens for $a \binom{b}{2,2}$ choices, the expression (\ref{uuuu}) 
in this case vanishes as $\varepsilon^{1+1/b}$ ;

2) some two of points $Z_1, \dots, Z_4$, belonging to different pairs, lie in one collection of $b$ critical points, and two other points in some other such collection.  This happens for $4\binom{a}{2}\binom{b}{2}^2$  choices. In a majority of these cases, namely when the two segments connecting the points inside these pairs are not the diagonals of a parallelogram, the exponent vanishes also as $\varepsilon^{1+1/b}$.

3) this is the exceptional subclass of the previous case, when we have such a parallelogram. This happens for $2\binom{a}{2} \binom{b}{2}$ choices. The expression (\ref{uuuu}) is then equal to zero, which means in particular that the line in $\C^k$ spanned by a linear function $\varphi$ never is generic with respect to the pure Stokes' set of our deformation of $f$. Let us study what will happen when we add a generic function of degree $\geq 2$ to $\varphi_0$.  

Adding the monomials $x^r$ and $y^r$ with arbitrary coefficients keeps our functions $f-\varepsilon \varphi$ split into the sums of two functions depending on $x$ and $y$ only. The set of critical values of such a function is a Minkovski sum of such sets of these functions, and therefore still 
contains parallelograms. So, the first possibility to move the corresponding expressions (\ref{uuuu}) from zero is to add the monomials proportional to $xy$. To avoid the details, consider the particular case $\varphi = (a+1)x+(b+1)y + xy$. Any parallelogram of critical values of the function $f-\varepsilon \varphi_0$, $\varepsilon >0$, in this case consists of its values at four points equal to $u_{1,2}\varepsilon^{1/a} +v_{1,2}\varepsilon^{1/b}$, where $u_{1,2}$ are some two different values of $1^{1/a}$, and $v_{1,2}$ some two values of $1^{1/b}$. Adding the monomial $xy$ to $\varphi_0$ we get the additions to the corresponding critical values, which are proportional, in the first approximation, to the values of the monomial $\varepsilon xy$ at these critical points. The expression (\ref{uuuu}) therefore changes from zero to asymptotically $ ((u_1v_1 + u_2v_2)-(u_1v_2+u_2v_1))\varepsilon^{1+1/a+1/b} =(u_1-u_2)(v_1-v_2)\varepsilon^{1+1/a+1/b}$. Adding the monomials of higher degrees divisible by $xy$ makes additions of higher orders in $\varepsilon$ to this expression and thus does not change its order of decrease. Finally, we get that any of our (deformed) parallelograms gives the contribution $1+1/a+1/b$ to the exponent of $\Omega(\varepsilon)$. \medskip    

Summing up all these exponents with their multiplicities, and then dividing the sum by 2, we get the formula (\ref{anscM2}). \hfill $\Box$

\begin{remark} \rm It is not very difficult to continue these calculations, considering also the pure Stokes' sets for functions $x^a+y^b$ with other parities of $a$ and $b$.
\end{remark}

\section{Appendix: motivations in algorithmic singularity theory}
\label{motivat}

The calculation of multiplicities considered above is needed for the optimization of a program enumerating all topoogically different pertubations of real function singularities. 

 For a description of this program see \S V.8 of the book \cite{APLT}, however
the web reference given there leads to an obsolete version of it (written in 1984 and described first in Chapter 5 of \cite{AVGL}). The actual versions of this program are available by \\
\verb"https://www.hse.ru/mirror/pubs/share/185895886" (for singularities of corank $\leq 2$) and \\
\verb"https://www.hse.ru/mirror/pubs/share/185895827"
 (for singularities of arbitrary ranks). The further versions of the program, which will use also the results of the present paper, will occur
at the bottom of the page \verb"https://www.hse.ru/en/org/persons/1297545#sci".

The main idea of the algorithm is as follows. Any Morse perturbation of a complicated singularity can be described in the terms of topological invariants related with the set of its critical values (including the imaginary ones), such as the order in $\R^1$ of their real parts, the Morse indices of corresponding critical points, and intersection matrices of related vanishing cycles. Such a collection of topological data is called a ``virtual morsification''.
The standard topological surgeries (such as collisions of critical points or critical values) can be modelled on the level of reconstructions of these collections. The algorithm starts from the data of a real morsification and applies to it all sequences of such admissible reconstructions. The virtual morsification of any real one surely occurs at some step of this algorithm. However in the case of sufficiently complicated singularities these admissible sequences are not necessarily finite, because they can include the changes of bases of vanishing cycles defined by the rotations of imaginary critical values around one another, and the occurring set of intersection matrices may be not finite. The results of this article give us  a priori estimates of the number of steps, after which any existing topological type of real perturbations will be attained, so that we can stop our algorithm. Namely, the multiplicity of the caustic (the Maxwell set, the pure Stokes' set) gives us an upper bound for the necessary number of collisions of critical points (respectively, of collisions of critical values not related with the collisions of corresponding critical points; the number of changes of the base of vanishing cycles related with the rotations of imaginary critical values). The study of the mixed Stokes' set helps us in the situation when two real critical values undergo the Morse surgery, go into the complex domain, travel there somehow and then come again to the real line
in some other place among the real critical values.

In the theory of real analytic function singularities, the mixed Stokes' set considered above is represented by (and usually coincides with the algebraic closure of) the set of parameters $\lambda \in \R^k$ such that the complexification of the function $f_\lambda$ has two complex conjugate critical values, whose real parts coincide with its critical value at some real critical point. Correspondingly, the real version of the pure Stokes' set of such a deformation is the closure of the set of parameters $\lambda$ such that the complexification of $f_\lambda$ has two pairs of conjugate critical values with equal real parts.

For the geometry of real Stokes' sets of singularities of small codimensions, see \cite{BH}, where also a link to the literature discussing the physical applications of this notion is given.  

The estimates of the sufficient numbers of virtual surgeries in our algorithm are justified by the following easy statement.  

\begin{proposition}
\label{triv}
If the local multiplicity of a real algebraic hypersurface $X \subset \R^k$ at the point $0 \in \R^k$ is equal to $d$, then any two points of its complement in a neighborhood of $0$ can be connected by a generic path intersecting this hypersurface at most $d$ times.
\end{proposition}

{\it Proof.} Consider a line in $\R^k$, which is generic with respect to our hypersurface among all lines through the origin in $\R^k$. There is a segment in this line whose unique intersection point with $X$ is its middle point $0 \in \R^k$. Choose two small balls in $\R^k$ centered at the endpoints of this segment and not intersecting the hypersurface $X$. Given two points in $\R^k \setminus X$ close to $0$, we can connect them by two paths, any of which is transversal to $X$ and consists of five parts: the first and the last of these parts connect our two points in $\R^k \setminus X$ to some two points very close to $0$; the second and the fourth ones are the segments parallel to the chosen line and ending in one of two chosen balls, and the middle part connects the obtained two points inside this ball. The sum of intersection numbers of these two paths with $X$ is at most $2d$, therefore at least one of these two numbers is not greater than $d$. \hfill $\Box$ 

\end{document}